\documentclass{amsart}
\usepackage{latexsym}

\usepackage{amsfonts,amssymb,amsmath}
\usepackage{amsthm}
\usepackage{graphicx,color}
\usepackage{color}
\usepackage{amsmath} 
\newtheorem{theorem}{Theorem}	
\newtheorem{lemma}{Lemma}[section]
\newtheorem{corollary}{Corollary}

\newtheorem{question}{Question}[section]	
\newfont{\bg}{cmr9 scaled\magstep2}
\newcommand{\bigzerol}{\smash{\lower1.0ex\hbox{\bg 0}}}

\makeatletter
\newcommand{\tpitchfork}{%
  \vbox{
    \baselineskip\z@skip
    \lineskip-.52ex
    \lineskiplimit\maxdimen
    \m@th
    \ialign{##\crcr\hidewidth\smash{$-$}\hidewidth\crcr$\pitchfork$\crcr}
  }%
}
\makeatother
\title[Geometric interpretation of generalized distance-squared mappings ]
{
Geometric interpretation of generalized distance-squared mappings 
of $\mathbb{R}^2$ into $\mathbb{R}^\ell$ $(\ell \geq 3)$
}
\author{Shunsuke Ichiki}
\thanks{The author is Research Fellow DC1 of Japan Society for the Promotion of Science }
\address{
Graduate School of Environment and Information Sciences,  
Yokohama National University, 
Yokohama 240-8501, Japan}
\email{ichiki-shunsuke-jb@ynu.jp}

\begin{document}
\date{}
%%%%%%%%%%%%%%%%%%%%%%%%%%%%%%%%%%%%%%%%%%%%% 
\begin{abstract}
%In \cite{G2}, generalized distance-squared mappings of the plane into the plane 
%having a generic central point are investigated. 
%Moreover, in the case that matrices $A$, which are 
%constructed by coefficients of generalized distance-squared mappings, are full, 
%geometric interpretation is also given. 
Generalized distance-squared mappings  
are quadratic mappings of $\mathbb{R}^m$ into $\mathbb{R}^\ell$ of special type. 
In the case that matrices $A$ 
constructed by coefficients of generalized distance-squared mappings 
of $\mathbb{R}^2$ into $\mathbb{R}^\ell$  ($\ell \geq3$) are full rank, 
the generalized distance-squared mappings having a generic central point have the following properties. 
In the case of $\ell=3$, they have only one singular point. 
On the other hand, in the case of $\ell>3$, they have no singular points. 
Hence, in this paper, the reason why in the case of $\ell=3$ 
(resp., in the case of $\ell>3$), 
they have only one singular point (resp., no singular points) is explained 
by giving a geometric interpretation to these phenomena.  
\end{abstract}
%%%%%%%%%%%%%%%%%%%%%%%%%%%%%%%%%%%%%%%%%%%%%%  
\maketitle
\noindent
\section{Introduction}\label{section 1}
%%%%%%%%%%%%%%%%%%%%%%%%%%%%%%%%%%%%%%%%%%%%%%%%%% 
%%%%%%%%%%%%%%%%%%%%%%%%%%%%%%%%%%%%%%%%%%%%%%%%%% 
Throughout this paper, positive integers are expressed by $i$, $j$, $\ell$ and $m$. 
In this paper, unless otherwise stated, all mappings belong to class $C^{\infty}$. 
Two mappings $f : \mathbb{R}^{m}\to \mathbb{R}^\ell$ and 
$g:\mathbb{R}^{m}\to \mathbb{R}^\ell$ are said  
to be  {\it  $\mathcal{A}$-equivalent} if 
there exist two diffeomorphisms $h : \mathbb{R}^{m}\to  \mathbb{R}^{m}$ and 
$H : \mathbb{R}^{\ell}\to  \mathbb{R}^{\ell}$ satisfying 
$f=H\circ g \circ h^{-1}$. 
Let $p_i=(p_{i1}, p_{i2}, \ldots, p_{im})$  $(1\le i\le \ell)$ 
(resp., $A=(a_{ij})_{1\le i\le \ell, 1\le j\le m}$) 
be a point of $\mathbb{R}^m$ 
(resp., an $\ell \times m$ matrix with non-zero entries).
Set $p=(p_1,p_2,\ldots,p_{\ell})\in (\mathbb{R}^m)^{\ell}$. 
Let $G_{(p, A)}:\mathbb{R}^m \to \mathbb{R}^\ell$ be the mapping 
defined by 
{\small 
\[
G_{(p, A)}(x)=\left(
\sum_{j=1}^m a_{1j}(x_j-p_{1j})^2, 
\sum_{j=1}^m a_{2j}(x_j-p_{2j})^2, 
\ldots, 
\sum_{j=1}^m a_{\ell j}(x_j-p_{\ell j})^2
\right), 
\] }where $x=(x_1, x_2, \ldots, x_m)\in \mathbb{R}^m$. 
The mapping $G_{(p, A)}$ is called a {\it generalized distance-squared mapping}, 
and the $\ell$-tuple of points $p=(p_1,\ldots ,p_{\ell})\in (\mathbb{R}^m)^{\ell}$ 
is called the {\it central point} 
of the generalized distance-squared mapping $G_{(p,A)}$. 
For a given matrix $A$, a property of 
generalized distance-squared mappings will be said to be 
true for a generalized distance-squared mapping 
having a {\it generic central point} 
if there exists a subset $\Sigma $ with Lebesgue measure zero 
of $(\mathbb{R}^m)^\ell$ such that for any $p \in (\mathbb{R}^m)^\ell-\Sigma$, 
the mapping $G_{(p,A)}$ has the property. 
%In this section, as an application of Theorem \ref{@f}, 
%we investigate the stability of the generalized distance-squared mappings. 
A {\it distance-squared mapping} $D_p$ 
(resp., {\it Lorentzian distance-squared mapping}
$L_p$) is the mapping $G_{(p,A)}$ 
satisfying that each entry of $A$ is $1$ 
(resp., $a_{i1}=-1$ and $a_{ij}=1$ $(j\ne 1)$). 
In \cite{D} (resp., \cite{L}), 
a classification result on distance-squared mappings $D_p$ 
(resp., Lorentzian distance-squared mappings
$L_p$) is given. Moreover, in \cite{G1} (resp., \cite{G2}), a classification result on generalized distance-squared mappings $G_{(p,A)}$ of $\mathbb{R}^2$ into  $\mathbb{R}^2$ 
(resp., $\mathbb{R}^{m+1}$ into  $\mathbb{R}^{\ell}$ ($\ell \geq 2m+1$)) is given. 

The important motivation for these investigations is as follows. 
Height functions and distance-squared functions have been investigated 
in detail so far, 
and they are well known as a useful tool 
in the applications of singularity theory to differential geometry 
(for example, see \cite{CS}). 
The mapping in which each component is a height function is 
nothing but a projection. 
Projections as well as height functions or distance-squared functions 
have been investigated so far. 

On the other hand, 
the mapping in which each component 
is a distance-squared function is a distance-squared mapping. 
Besides, the notion of generalized distance-squared mappings is 
an extension of the distance-squared mappings. 
Therefore, it is natural to investigate generalized distance-squared mappings 
as well as projections. 

In \cite{G1}, 
a classification result on 
generalized distance-squared mappings of $\mathbb{R}^2$ into $\mathbb{R}^2$ 
is given. 
If the rank of $A$ is two,  
the generalized distance-squared mapping 
having a generic central point 
is a mapping of which any singular point is 
a {\color{black}{\rm fold point}} 
except one {\color{black}{\rm cusp point}}.     
%(for details on fold points and cusp points, refer to \cite{plane to plane}). 
The singular set is a rectangular hyperbola. 
Moreover, in \cite{G1}, a 
geometric interpretation of a singular set of 
generalized distance-squared mappings of $\mathbb{R}^2$ into $\mathbb{R}^2$ having a generic central point is also given in the case of ${\rm rank\ }A=2$.  
By the geometric interpretation,  
the reason why the mappings have only one cusp point is explained.  

%If the rank of $A$ is one, 
%a generalized distance-squared mapping 
%having a generic central point is $\mathcal{A}$-equivalent 
%to 
%{\color{black}{\rm the normal form of fold singularity}} $(x_1,x_2)
%\mapsto (x_1,x_2^2)$. 

On the other hand, in \cite{G2}, 
a classification result on 
generalized distance-squared mappings 
of $\mathbb{R}^{m+1}$ into $\mathbb{R}^{\ell}$ $(\ell \geq 2m+1)$ 
is given. 
%If the rank of $A$ is $m+1$, 
%a generalized distance-squared mapping 
%having a generic central point 
%is $\mathcal{A}$-equivalent to 
%{\color{black}{\rm the normal form of Whitney umbrella}} 
%%as follows:
%$
%(x_1,\ldots ,x_{m+1})\mapsto (x_1^2,x_1x_2,\ldots ,x_1x_{m+1},x_2,\ldots ,x_{m+1})
%$.
%If the rank of $A$ is less than $m+1$, 
%a generalized distance-squared mapping of $\mathbb{R}^{m+1}$ into 
%$\mathbb{R}^{2m+1}$ 
%having a generic central point 
%is $\mathcal{A}$-equivalent to the inclusion 
%%as follows:
%$
%(x_1,\ldots ,x_{m+1})\mapsto (x_1,\ldots ,x_{m+1},0,\ldots ,0)
%$. 
%\par
As the special case of $m=1$, we have the following. 
\begin{theorem}[\cite{G2}]\label{main}
Let $\ell$ be an integer satisfying $\ell \geq 3$. 
Let $A=(a_{ij})_{1\leq i \leq \ell, 1\leq j \leq 2}$ be an $\ell \times 2$ matrix with 
non-zero entries satisfying ${\rm rank\ }A=2$. 
Then, the following hold: 
\begin{enumerate}
\item 
In the case of $\ell=3$, 
there exists a proper algebraic subset $\Sigma _A\subset (\mathbb{R}^2)^3$ 
such that for any $p\in  (\mathbb{R}^2)^3-\Sigma _A$, 
the mapping $G_{(p,A)}$ is $\mathcal{A}$-equivalent to 
the normal form of Whitney umbrella
$(x_1,x_2)\mapsto (x_1^2, x_1x_2, x_2)$. 
\item 
In the case of $\ell>3$, 
there exists a proper algebraic subset $\Sigma _A\subset (\mathbb{R}^2)^\ell$ 
such that for any $p\in  (\mathbb{R}^2)^\ell-\Sigma _A$, 
the mapping $G_{(p,A)}$ is $\mathcal{A}$-equivalent to the inclusion 
$(x_1,x_2)\mapsto (x_1, x_2, 0,\ldots ,0)$. 
\end{enumerate}
\end{theorem}

As described above, in \cite{G1}, a geometric interpretation of 
generalized distance-squared mappings of $\mathbb{R}^2$ into $\mathbb{R}^2$ 
having a generic central point is given in the case that the matrix $A$ is full rank.   
On the other hand, in this paper, 
a geometric interpretation of 
generalized distance-squared mappings of $\mathbb{R}^2$ into $\mathbb{R}^\ell$ 
$(\ell \geq3)$  
having a generic central point is given in the case that the matrix $A$ is full rank 
(for the reason why we concentrate on the case that the matrix $A$ is full rank, 
see Remark \ref{not full rank}). 
Hence, by \cite{G1} and this paper, 
geometric interpretations of generalized distance-squared mappings of the plane 
having a generic central point in the case that the matrix $A$ is full rank 
are completed. 

The main purpose of this paper is to give a geometric interpretation of  
Theorem \ref{main}. 
Namely, the main purpose of this paper is to answer the following question. 
\begin{question}\label{question 1}
Let $\ell$ be an integer satisfying $\ell \geq 3$. 
Let $A=(a_{ij})_{1\leq i \leq \ell, 1\leq j \leq 2}$ be an $\ell \times 2$ matrix with 
non-zero entries satisfying ${\rm rank\ }A=2$. 
\begin{enumerate}
\item 
In the case of $\ell=3$, 
why do generalized distance-squared mappings 
$G_{(p,A)}:\mathbb{R}^2\to \mathbb{R}^3$ having a generic central point have only one singular point ?
\item 
On the other hand, 
in the case of $\ell>3$, 
why do generalized distance-squared mappings 
$G_{(p,A)}:\mathbb{R}^2\to \mathbb{R}^\ell$ having a generic central point have 
no singular points ?
\end{enumerate}
\end{question}
%
%We have another original motivation. 
%Height functions and distance-squared functions have been investigated 
%in detail {\color{black}so far}, 
%and they are a useful tool 
%in the applications of singularity theory to differential geometry 
%({\color{black}for instance}, see \cite{CS}). 
%The mapping in which each component is a height function is 
%{\color{black}nothing but} a projection. 
%In \cite{GP}, compositions of 
%{\color{black}generic} projections and embeddings are investigated.
%
%On the other hand, 
%the mapping in which each component 
%is a distance-squared function is a distance-squared mapping. 
%{\color{black}And}, the notion of generalized distance-squared mapping is 
%an extension of the distance-squared mappings. 
%Therefore, 
%{\color{black}it is again natural to} investigate compositions 
%{\color{black}with generic} generalized distance-squared mappings. 
\subsection{Remark}\label{not full rank}
In the case that the matrix $A$ is not full rank $({\rm rank\ }A=1)$, 
for any $\ell\geq 3$, 
the generalized distance-squared mapping 
of  $\mathbb{R}^2$ into $\mathbb{R}^\ell$ having a generic central point  
is $\mathcal{A}$-equivalent to only the inclusion $(x_1,x_2)\mapsto (x_1, x_2, 0, \ldots , 0)$ (see Theorem 3 in \cite{G2}). 
On the other hand, in the case that the matrix $A$ is full rank $({\rm rank\ }A=2)$, 
the phenomenon in the case of $\ell=3$ is completely different from the phenomenon  in 
the case of $\ell>3$ (see Theorem \ref{main}). 
Hence, we concentrate on the case that the matrix $A$ is full rank . 
\par 
\bigskip 
%%%%%%%%%%%  
In Section \ref{section 2}, 
some assertions and definitions are prepared for 
answering Question \ref{question 1}, and 
the answer to Question \ref{question 1} is stated. 
In Section \ref{section 4}, the proof of a lemma of Section \ref{section 2} is given.    
%%%%%%%%%%%%%%%%%%%%%%%%%%%%%%%%%%%%%%%%%%%%%%%%%%
%%%%%%%%%%%%%%%%%%%%%%%%%%%%%%%%%%%%%%%%%%%%%%%%%%    
\section{The answer to Question \ref{question 1}}\label{section 2}
%%%%%%%%%%%%%%%%%%%%%%%%%%%%%%%%%%%%%%%%%%%%%%%%%%
%%%%%%%%%%%%%%%%%%%%%%%%%%%%%%%%%%%%%%%%%%%%%%%%%% 
Firstly, in order to answer Question \ref{question 1}, 
some assertions and definitions 
 are prepared. 
By Theorem \ref{main}, it is clearly seen that the following assertion holds. 
The assertion is important for giving an geometric interpretation. 
\begin{corollary}\label{isomorphism}
Let $\ell$ be an integer satisfying $\ell \geq 3$. 
Let $A=(a_{ij})_{1\leq i \leq \ell, 1\leq j \leq 2}$ 
$(resp., B=(b_{ij})_{1\leq i \leq \ell, 1\leq j \leq 2})$ be an $\ell \times 2$ matrix with 
non-zero entries 
satisfying {\rm rank\ }$A=2$ $($resp., {\rm rank\ }$B=2)$.  
Then, there exist  proper algebraic subsets $\Sigma_A$ and $\Sigma_B$ of 
$(\mathbb{R}^2)^\ell$ such that for any $p \in (\mathbb{R}^2)^\ell-\Sigma_A$ 
and for any $q \in (\mathbb{R}^2)^\ell-\Sigma_B$, the mapping 
$G_{(p,A)}$ is $\mathcal{A}$-equivalent to the mapping $G_{(q, B)}$. 
\end{corollary}
Let $F_p:\mathbb{R}^2\to \mathbb{R}^\ell$ $(\ell \geq3)$ be the mapping 
defined by 
\begin{eqnarray*}
&&F_p(x_1,x_2)\\
&=&
\left(
a(x_1-p_{11})^2+b(x_2-p_{12})^2, (x_1-p_{21})^2+(x_2-p_{22})^2, \ldots , 
\right. 
\\
&&
\left. 
(x_1-p_{\ell 1})^2+(x_2-p_{\ell 2})^2
\right),
\end{eqnarray*}
where $0<a<b$ and $p=(p_{11},p_{12}, \ldots ,p_{\ell 1},p_{\ell 2})$.  
Remark that the mapping $F_p$ is the generalized distance-squared mapping $G_{(p,B)}$,  
where 
\begin{eqnarray*}
B=
\left(
    \begin{array}{ccc}
      a& b\\
      1& 1 \\
       \vdots & \vdots \\
      1& 1 \\
    \end{array}
\right). 
\end{eqnarray*}
 
Since the rank of the matrix $B$ is two, 
by Corollary \ref{isomorphism}, in order to answer 
Question \ref{question 1}, it is sufficient to answer the following question. 

\begin{question}\label{question 2}
Let $\ell$ be an integer satisfying $\ell \geq 3$. 
\begin{enumerate}
\item 
In the case of $\ell=3$, 
why do the mappings $F_p : \mathbb{R}^2\to \mathbb{R}^3$ having a generic central point have only one singular point ?
\item 
On the other hand, 
in the case of $\ell>3$, 
why do the mappings 
$F_p : \mathbb{R}^2\to \mathbb{R}^\ell$ having a generic central point have 
no singular points ?
\end{enumerate}
\end{question}

The mapping $F_p=(F_1,\ldots ,F_\ell)$ determines $\ell$-foliations 
$\mathcal{C}_{p_1}(c_1), \ldots ,\mathcal{C}_{p_\ell}(c_\ell)$ 
in the plane defined by 
\begin{eqnarray*}
\mathcal{C}_{p_i}(c_i)=\{(x_1, x_2)\in \mathbb{R}^2 \mid F_i(x_1,x_2)=c_i\}, 
%\mathcal{C}_{p_2}(c_2)=\{(x,y)\mid (x_1-p_{21})^2+(x_2-p_{22})^2=c_2\}, \\
%\mathcal{C}_{p_3}(c_1)=\{(x,y)\mid (x_1-p_{31})^2+(x_2-p_{32})^2=c_3\}, 
\end{eqnarray*}
where $c_i\geq 0$ $(1\leq i \leq \ell)$ and 
$p=(p_1,\ldots ,p_\ell)\in (\mathbb{R}^2)^\ell$. 
%$L$-foliations are tangent at a point if and only if the point is a singular point of the mapping $F_p$. 
For a given central point $p=(p_1,\ldots ,p_\ell)\in (\mathbb{R}^2)^\ell$, 
a point $q\in \mathbb{R}^2$ is a singular point of the mapping $F_p$ 
if and only if the $\ell$-foliations $\mathcal{C}_{p_i}(c_i)$ $(1\leq i \leq \ell)$ 
defined by the point $p$ are tangent at the point $q$, 
where $(c_1,\ldots ,c_\ell)=F_p(q)$. 

For a given central point $p=(p_1,\ldots ,p_\ell)\in (\mathbb{R}^2)^\ell$, in the case that a point $q\in \mathbb{R}^2$ is a singular point of 
the mapping $F_p$, there may exist an integer $i$ such that the foliation $\mathcal{C}_{p_i}(c_i)$ is merely a point, where $c_i=F_i(q)$. 
However, by the following lemma, 
we see that the trivial phenomenon seldom occurs 
(for the proof of Lemma \ref{omit}, see Section \ref{section 4}). 

%The following proposition gives the relation of a singular point of the mappings $F_p$ 
%and tangent points of 
%$\ell$-foliations $\mathcal{C}_{p_1}(c_1), \ldots, \mathcal{C}_{p_\ell}(c_\ell)$ 
%(for the proof of Proposition \ref{singular tangent}, see Section \ref{section 4}). 
%%\begin{proposition}\label{singular tangent}
%A point $(x_1,x_2)$ is a singular point of the mapping $F_p$ and $(x_1,x_2)\not=p_1,\ldots, p_\ell$ 
%if and only if $\ell$-foliations 
%%either $(1)$ or $(2)$ holds:
%%\begin{itemize}
%%\item[$(1)$] 
%$\mathcal{C}_{p_1}(c_1),\ldots ,\mathcal{C}_{p_\ell}(c_\ell)$ are tangent at the point $(x_1,x_2)$.
%%\item[$(2)$] 
%%the point $(x,y)$ is $p_1$, $p_2$ or $p_3$ and the following holds: 
%%\begin{eqnarray*}
%%  {\rm rank} \left(
%%    \begin{array}{rrr}
%%      a(x-p_{11}) & b(y-p_{12}) \\
%%      x-p_{21} & y-p_{22} \\
%%      x-p_{31} & y-p_{32}
%%    \end{array}
%%  \right)<2
%%\end{eqnarray*}
%%\end{itemize}
%
%\end{proposition}

%We also get the following lemma (for the proof of Lemma \ref{omit}, see 
%Section \ref{section 5}). 
\begin{lemma}\label{omit}
Let $\ell$ be an integer satisfying $\ell \geq 3$. 
Then, there exists a proper algebraic subset $\Sigma \subset (\mathbb{R}^2)^\ell$ such that 
for any central point $p=(p_1,\ldots ,p_\ell)\in (\mathbb{R}^2)^\ell-\Sigma$, 
if a point $q\in \mathbb{R}^2$ is a singular point of the mapping $F_p$,  then 
the $\ell$-foliations $\mathcal{C}_{p_1}(c_1)$ and $\mathcal{C}_{p_i}(c_i)$ $(2\leq i \leq \ell)$ are an ellipse and $(\ell-1)$-circles, respectively, where $(c_1,\ldots ,c_\ell)=F_p(q)$.   
\end{lemma}
%By combining Proposition \ref{singular tangent} and Lemma \ref{omit}, 
%we have the following. 
%Remark that the subset $\Sigma$ in Theorem \ref{main2} is the subset $\Sigma$ 
%in Lemma \ref{omit}.  
%\begin{theorem}\label{main2}
%There exists a proper algebraic subset $\Sigma \subset (\mathbb{R}^2)^\ell$ 
%such that for any $p\in (\mathbb{R}^2)^\ell-\Sigma $, 
%the following assertions are equivalent. 
%\begin{itemize}
%\item[$(1)$] 
%a point $(x_1,x_2)$ is a singular point of the mapping $F_p$. 
%\item[$(2)$] 
%$\ell$-foliations $\mathcal{C}_{p_1}(c_1), \ldots ,\mathcal{C}_{p_\ell}(c_\ell)$ are tangent at $(x_1,x_2)$.
%\end{itemize}
%\end{theorem}
\subsection{Answer to Question \ref{question 1}}\label{answer}
As described above, in order to answer Question \ref{question 1}, 
it is sufficient to answer Question \ref{question 2}. 
\begin{enumerate}
\item 
We will answer (1) of Question \ref{question 2}. 
The phenomenon that the mapping $F_p:\mathbb{R}^2\to \mathbb{R}^3$ having 
a generic central point has only one singular point  
can be explained by the following geometric interpretation. Namely,  
constants $c_i \geq0$ $(i=1,2,3)$ such that three foliations 
 $\mathcal{C}_{p_1}(c_1), \mathcal{C}_{p_2}(c_2)$  and $\mathcal{C}_{p_3}(c_3)$ defined by the central point $p=(p_1, p_2, p_3)\in (\mathbb{R}^2)^3$ are tangent are uniquely determined, 
 and the tangent point is also unique. 
Moreover, in the case, remark that by Lemma \ref{omit}, 
the three foliations $\mathcal{C}_{p_1}(c_1), \mathcal{C}_{p_2}(c_2)$  and $\mathcal{C}_{p_3}(c_3)$ defined by almost all (in the sense of Lebesgue measure) $(p_1,p_2,p_3)\in (\mathbb{R}^2)^3$ are an ellipse and two circles, respectively.  

Furthermore, by the geometric interpretation, 
we can also see the location of the singular point of the mapping $F_p$ 
having a generic central point
(for example, see Figure \ref{figure 1}). 
\item 
We will answer (2) of Question \ref{question 2}. 
The phenomenon that the mapping $F_p:\mathbb{R}^2\to \mathbb{R}^\ell$ $(\ell >3)$ having 
a generic central point has no singular points 
can be explained by the following geometric interpretation. Namely, 
for any constants $c_i\geq0$ $(1\leq i \leq \ell)$, 
$\ell$-foliations $\mathcal{C}_{p_1}(c_1), \ldots , \mathcal{C}_{p_\ell}(c_\ell)$ defined by the central point $p=(p_1, \ldots , p_\ell)\in (\mathbb{R}^2)^\ell$ are not tangent at any points 
(for example, see Figure \ref{figure 2}).  
%Moreover, in the case, remark that by Lemma \ref{omit}, 
%the $\ell$-foliations $\mathcal{C}_{p_1}(c_1)$ and $\mathcal{C}_{p_i}(c_i)$ $(2\leq i \leq \ell)$ defined by almost all (in the sense of Lebesgue measure) $p\in (\mathbb{R}^2)^\ell$ are an ellipse and circles, respectively.  

\begin{figure}[htbp]
\begin{center}
\includegraphics[width=.68\linewidth]{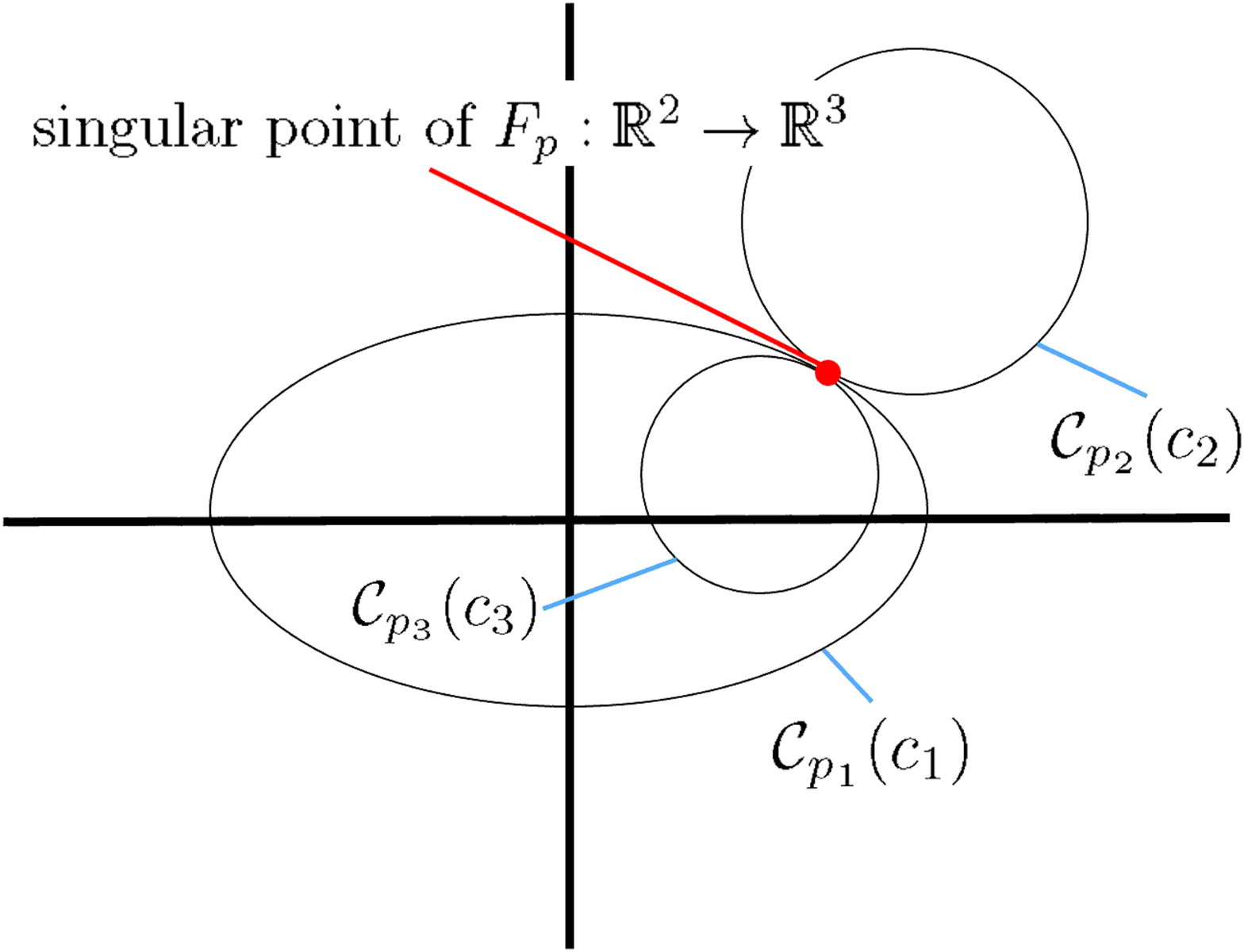}
\caption{A figure for (1) of Subsection \ref{answer}}
\label{figure 1}
\end{center}
\end{figure}  
\newpage
\begin{figure}[htbp]
\begin{center}
\includegraphics[width=.68\linewidth]{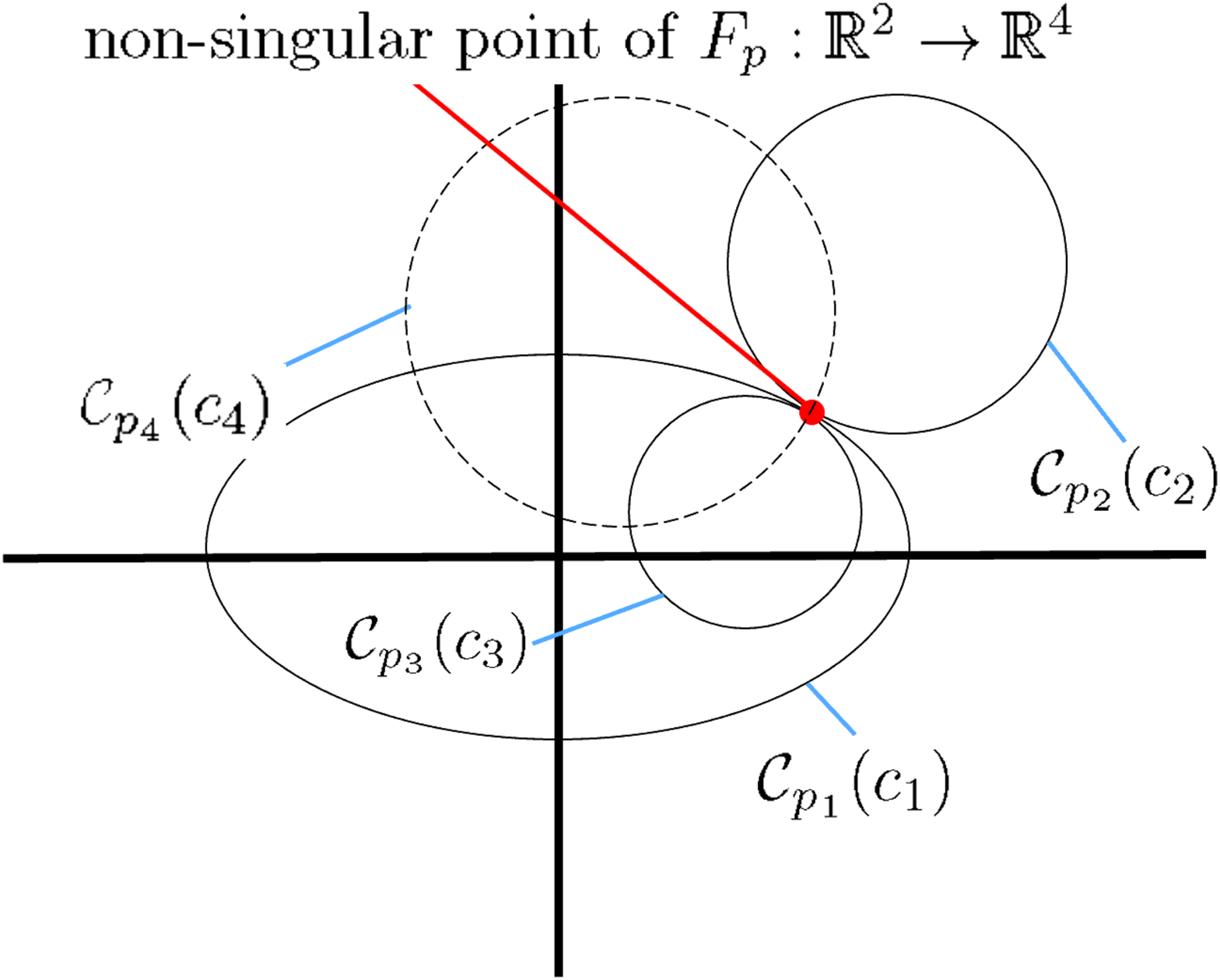}
\caption{A figure for (2) of Subsection \ref{answer}}
\label{figure 2}
\end{center}
\end{figure}  
\end{enumerate}

\subsection{Remark}
The geometric interpretation (the answer to Question \ref{question 1}) has one more advantage. 
By the interpretation, we get the following assertion from the viewpoint of the contacts amongst one ellipse and some circles. 
\begin{corollary}
Let $a, b$ be real numbers satisfying $0<a<b$.  
\begin{enumerate}
\item 
There exists a proper algebraic subset $\Sigma$ of $(\mathbb{R}^2)^3$ 
such that for any $(p_1,p_2,p_3)\in (\mathbb{R}^2)^3-\Sigma$, 
constants $c_i>0$ $(i=1,2,3)$ such that one ellipse $a(x_1-p_{11})^2+b(x_2-p_{12})^2=c_1$ and two circles $(x_1-p_{i1})^2+(x_2-p_{i2})^2=c_i$ $(i=2, 3)$ are 
tangent are uniquely determined, where $p_i=(p_{i1}, p_{i2})$. Moreover, the tangent point 
is also unique. 
\item 
On the other hand, in the case of $\ell>3$, 
there exists a proper algebraic subset $\Sigma$ of $(\mathbb{R}^2)^\ell$ 
such that for any $(p_1,\ldots, p_\ell)\in (\mathbb{R}^2)^\ell-\Sigma$, 
for any $c_i>0$ $(i=1,\ldots ,\ell)$, 
 the one ellipse $a(x_1-p_{11})^2+b(x_2-p_{12})^2=c_1$ and the $(\ell-1)$-circles   
 $(x_1-p_{i1})^2+(x_2-p_{i2})^2=c_i$ $(i=2,\ldots ,\ell)$ are not tangent 
 at any points, where $p_i=(p_{i1}, p_{i2})$.  
\end{enumerate}
\end{corollary}
\section{Proof of Lemma \ref{omit}} \label{section 4}
%%%%%%%%%%%%%%%%%%%%%%%%%%%%%%%%%%%%%%%%%%%%%%%%%%
%%%%%%%%%%%%%%%%%%%%%%%%%%%%%%%%%%%%%%%%%%%%%%%%%% 
The Jacobian matrix of the mapping $F_p$ at $(x_1,x_2)$ is the following.  
\begin{eqnarray*}
JF_{p{(x_1,x_2)}} =2\left(
    \begin{array}{ccc}
 a(x_1-p_{11}) & b(x_2-p_{12}) \\
      \vdots & \vdots \\
      x_1-p_{\ell 1} & x_2-p_{\ell 2}
    \end{array}
  \right). 
\end{eqnarray*}
Let $\Sigma_i$ be a subset of $(\mathbb{R}^2)^\ell$ consisting of 
$p=(p_1,\ldots ,p_\ell)\in (\mathbb{R}^2)^\ell$ satisfying    
$p_i\in \mathbb{R}^2$ is a singular point of $F_p$ $(1\leq i \leq \ell)$. 
Namely, 
for example, $\Sigma_1$ is the subset of $(\mathbb{R}^2)^\ell$ consisting of 
$p=(p_1,\ldots ,p_\ell)\in (\mathbb{R}^2)^\ell$ satisfying 
\begin{eqnarray*}
{\rm rank\ } \left(
    \begin{array}{cccc}
    0 & 0 \\
 p_{11}-p_{21} & p_{12}-p_{22} \\
      \vdots & \vdots \\
      p_{11}-p_{\ell 1} & p_{12}-p_{\ell 2}
    \end{array}
  \right)<2. 
\end{eqnarray*}
By $\ell \geq 3$, it is clearly seen that $\Sigma_1$ is a proper algebraic subset of 
$(\mathbb{R}^2)^\ell$. 
Similarly, for any $i$ $(2\leq i \leq \ell)$, we see that  $\Sigma_i$ is also a proper algebraic subset 
of 
$(\mathbb{R}^2)^\ell$. 
Set $\Sigma=\cup_{i=1}^\ell \Sigma_i$. Then, $\Sigma$ is also a proper algebraic subset 
of $(\mathbb{R}^2)^\ell$. 

Let $p=(p_1,\ldots ,p_\ell)\in (\mathbb{R}^2)^\ell-\Sigma$ be a central point, and let $q$ be 
a singular point of the mapping $F_p$ defined by the central point. 
Then, suppose that there exists an integer $i$ such that 
the foliation $\mathcal{C}_{p_i}(c_i)$ is not an ellipse or a circle, 
where $c_i=F_i(q)$ $(F_p=(F_1, \ldots ,F_\ell))$. 
Then, we get $c_i=0$. Hence, we have $q=p_i$. 
This contradicts the assumption $p\in (\mathbb{R}^2)^\ell-\Sigma$.

%
%%$\Sigma_i=\{p\in (\mathbb{R}^2)^3\mid {\rm rank} JF_{p{(p_{i1},p_{i2})}} <2 \}$
%We will show that $\Sigma_i$ is a proper algebraic set $(i=1,2,3)$.  

\hfill $\Box$
\section*{Acknowledgements}
The author is grateful to Takashi Nishimura for his kind advices. 
The author is supported by JSPS KAKENHI Grant Number 16J06911.
%%%%%%%%%%%%%%%%%%%%%%%%%%%%%%%%%%%%%%%%%%%%%%%%%%   
%%%%%%%%%%%%%%%%%%%%%%%%%%%%%%%%%%%%%%%%%%%%%%%%%%  

\end{document}